\documentclass[11pt]{article}
\usepackage{amsmath,amsthm}

\def\refeq#1{\if\workingver y(\ref{#1})-[[#1]]\else(\ref{#1})\fi}
\def\refth#1{\if\workingver y\ref{#1}-[[#1]]\else\ref{#1}\fi}
\def\mylabel#1{\if\workingver y\label{#1}{\bf\ \ [[#1]]\ \ }
\else\label{#1}\fi}
\def\mybibitem#1{\if\workingver y\bibitem{#1}{\bf\ \ [[#1]]\ \ }
\else\bibitem{#1}\fi}

\def\institute#1{\gdef\@institute{#1}}

\def\institutename{\par
 \begingroup
 \parskip=\z@
 \parindent=\z@
 \setcounter{@inst}{1}%
 \def\and{\par\stepcounter{@inst}%
 \noindent$^{\the@inst}$\enspace\ignorespaces}%
 \setbox0=\vbox{\def\thanks##1{}\@institute}%
 \ifnum\c@@inst=1\relax
 \else
   \setcounter{footnote}{\c@@inst}%
   \setcounter{@inst}{1}%
   \noindent$^{\the@inst}$\enspace
 \fi
 \ignorespaces
 \@institute\par
 \endgroup}

\newtheorem{thm}{Theorem}
\newtheorem{lem}[thm]{Lemma}

\newtheorem{cor}[thm]{Corollary}

\newtheorem{rem}[thm]{Remark}

\def\begeq#1{\begin{equation}\mylabel{#1}}
\def\endeq{\end{equation}}

\def\begalg{\begin{alg}}
\def\endalg{\end{alg}}

\let\workingver=n

\newcommand{\beq}{\begin{equation}}

\def\sms{\small\scshape}

\def\al{\alpha}
\def\be{\beta}

\begin{document}

\title{Generating Functions, Weighted and Non-Weighted
Sums for Powers of Second-Order Recurrence Sequences}

\author{Pantelimon St\u anic\u a \vspace{.2cm}
\thanks{\em
Also associated with the Institute of
Mathematics of Romanian Academy,
Bucharest, Romania}\\
\small Auburn University Montgomery,
\small Department of Mathematics\\
\small Montgomery, AL 36117, USA\\
\small e-mail: \em  stanpan@strudel.aum.edu}
\date{
}
\maketitle

\begin{abstract}
\baselineskip=2\baselineskip
In this paper we find closed forms of the generating function
$\displaystyle \sum_{k=0}^\infty U_n^r x^n$, for powers of any
non-degenerate second-order recurrence
sequence,  $U_{n+1}=aU_{n}+bU_{n-1},\ a^2+4b\not= 0$, completing a study
began by Carlitz \cite{Carlitz1962} and Riordan \cite{Riordan1962} in 1962.
Moreover, we generalize a theorem of Horadam \cite{Horadam1994} on
 partial sums involving
such sequences. Also, we find closed forms for weighted (by binomial
coefficients) partial sums
of powers of any non-degenerate second-order recurrence sequences.
 As corollaries we give some known and seemingly unknown
identities and derive some very interesting congruence relations
involving Fibonacci and Lucas sequences.
\end{abstract}

\section{Introduction}

\baselineskip=2\baselineskip

DeMoivre (1718) used the generating function (found by using the
recurrence) for the Fibonacci sequence \(\displaystyle
\sum_{i=0}^\infty F_i x^i=\frac{x}{1-x-x^2} \), to obtain the
identities
$F_n=\frac{\al^n-\be^n}{\sqrt{5}}, L_n=\al^n+\be^n$ ({\em Lucas
numbers}) with  $\al=\frac{1+\sqrt{5}}{2},
\be=\frac{1-\sqrt{5}}{2}$, called {\em Binet formulas},
in honor of Binet who in fact rediscovered them more than one
hundred years later, in 1843 (see \cite{Vajda1989}).
 Reciprocally, using the Binet
formulas, we can find the generating function easily
\(\displaystyle \sum_{i=0}^\infty F_i x^i =\frac{1}{\sqrt{5}}
\sum_{i=0}^\infty (\al^i-\be^i)x^i =\frac{1}{\sqrt{5}}\left(
\frac{1}{1-\al x}-\frac{1}{1-\be x} \right)= \frac{x}{1-x-x^2}, \)
since $\al\be=-1, \al+\be=1$.

The question that arises is whether we can find a closed form for
the generating function for powers of Fibonacci numbers, or
better yet, for powers of any second-order recurrence sequences.
Carlitz \cite{Carlitz1962} and
 Riordan \cite{Riordan1962} were unable to find
 the closed form for the generating functions $F(r,x)$ of $F_n^r$,
 but found a recurrence relation among them, namely
 \[
\displaystyle
(1-L_r x+(-1)^r x^2) F(r,x)=1+rx\sum_{j=1}^{[\frac{r}{2}]}
(-1)^j \frac{A_{rj}}{j} F(r-2j,(-1)^jx),
 \]
 with $A_{rj}$ having a complicated structure (see also \cite{Horadam1965}).
 We are
 able to complete the study  began by them and
 find a closed form for the generating function for
 powers of any non-degenerate second-order recurrence sequence.
 We would like to point out,  that this
  "forgotten" technique we employ
 can be used to attack successfully other sums or series involving
  any second-order recurrence sequence.
  In this paper we also find  closed forms for non-weighted
  partial sums for non-degenerate second-order
  recurrence sequences, generalizing a theorem of Horadam
  \cite{Horadam1994} and also weighted (by the binomial coefficients)
 partial sums for such sequences.

\section{Generating Functions}

We consider the general non-degenerate second-order recurrences
$U_{n+1}=aU_n+bU_{n-1}$, $a^2+4b\not= 0$.
We intend to find the generating function
\(\displaystyle
U(r,x)=\sum_{i=0}^\infty U_i^r x^i.
\)
It is known that the Binet formula for the sequence $U_n$ is
$U_n=A\al^n-B\be^n$, where $\al=\frac{1}{2}(a+\sqrt{a^2+4b}),
\be=\frac{1}{2}(a-\sqrt{a^2+4b})$ and $A=\frac{U_1-U_0\be}{\al-\be},
B=\frac{U_1-U_0\al}{\al-\be}$.
We associate the sequence
$V_n=\al^n+\be^n$, which satisfies the same recurrence, with the initial
conditions $V_0=2,V_1=a$.
\begin{thm}
\label{thmU}
We have
\[
U(r,x)=\sum_{k=0}^{\frac{r-1}{2}} (-1)^k A^k B^k\binom{r}{k}
\frac{A^{r-2k}-B^{r-2k}+(-b)^k(B^{r-2k}\al^{r-2k}-
A^{r-2k} \be^{r-2k})x}{1- (-b)^k V_{r-2k} -x^2},
\]
if $r$ odd, and
\[
\begin{split}
U(r,x)=&
\sum_{k=0}^{\frac{r}{2}-1}(-1)^k A^k B^k \binom{r}{k}
\frac{B^{r-2k}+A^{r-2k}-(-b)^k(B^{r-2k}\al^{r-2k}+
A^{r-2k}\be^{r-2k})x}{1- (-b)^kV_{r-2k}x+x^2}\\
&+\binom{r}{\frac{r}{2}}\frac{A^{\frac{r}{2}}
(-B)^{\frac{r}{2}}}{1-(-1)^{\frac{r}{2}} x},\ \text{if $r$ even.}
\end{split}
\]
\end{thm}
\begin{proof}
We evaluate
\begin{eqnarray*}
U(r,x)
&=&
\sum_{i=0}^\infty \left( \sum_{k=0}^r \binom{r}{k}
(A\al^i)^k(-B\be^i)^{r-k} \right)x^i\\
&=&
\sum_{k=0}^r \binom{r}{k}A^k(-B)^{r-k}\sum_{i=0}^\infty
(\al^k \be^{r-k}x)^i
\\
&=&
\sum_{k=0}^r \binom{r}{k} A^k(-B)^{r-k} \frac{1}{1-\al^k \be^{r-k}x}.
\end{eqnarray*}
If $r$ odd, then associating $k\leftrightarrow r-k$, we get
\[
\begin{split}
U(r,x)
=&
\sum_{k=0}^{\frac{r-1}{2}} \binom{r}{k}
\left(\frac{A^k (-B)^{r-k}}{1-\al^k\be^{r-k}x} +
\frac{A^{r-k} (-B)^k}{1-\al^{r-k}\be^k x}\right)\\
=&
\sum_{k=0}^{\frac{r-1}{2}} (-1)^k \binom{r}{k}
\left(\frac{A^{r-k} B^k}{1-\al^{r-k}\be^k x}-
\frac{A^k B^{r-k}}{1-\al^k\be^{r-k}x}
\right)\\
=&
\sum_{k=0}^{\frac{r-1}{2}} (-1)^k \binom{r}{k}
\frac{A^{r-k}B^k-A^k B^{r-k}+(A^k B^{r-k}\al^{r-k}\be^k-
A^{r-k} B^k\al^k\be^{r-k})x}{1-(\al^k\be^{r-k}+\al^{r-k}\be^k)x+
\al^r\be^r x^2}\\
=&
\sum_{k=0}^{\frac{r-1}{2}} (-1)^k\binom{r}{k}
\frac{A^{r-k}B^k-A^k B^{r-k}+(-1)^kb^k(A^k B^{r-k}\al^{r-2k}-
A^{r-k} B^k\be^{r-2k})x}{1- (-1)^k b^k V_{r-2k} -x^2}.
\end{split}
\]
If $r$ even, then then associating $k\leftrightarrow r-k$,
except for the middle term, we get
\[
\begin{split}
U(r,x)
=&
\sum_{k=0}^{\frac{r}{2}-1} \binom{r}{k}
\left(\frac{A^k (-B)^{r-k}}{1-\al^k \be^{r-k}x} +
\frac{A^{r-k} (-B)^k}{1-\al^{r-k}\be^k x}\right)+
\binom{r}{\frac{r}{2}}\frac{A^{\frac{r}{2}}
(-B)^{\frac{r}{2}}}{1-(-1)^{\frac{r}{2}} x}\\
=&
\sum_{k=0}^{\frac{r}{2}-1} (-1)^k \binom{r}{k}
\left(\frac{A^{r-k} B^k}{1-\al^{r-k}\be^k x}+
\frac{A^k B^{r-k}}{1-\al^k\be^{r-k}x}
\right)
+ \binom{r}{\frac{r}{2}}\frac{A^{\frac{r}{2}}
(-B)^{\frac{r}{2}}}{1-(-1)^{\frac{r}{2}} x}
\end{split}
\]
\[
\begin{split}
=&
\sum_{k=0}^{\frac{r}{2}-1}(-1)^k \binom{r}{k}
\frac{A^k B^{r-k}+A^{r-k}B^k-(A^k B^{r-k}\al^{r-k}\be^k+
A^{r-k} B^k\al^k\be^{r-k})x}{1-(\al^k\be^{r-k}+\al^{r-k}\be^k)x+
\al^r\be^r x^2}\\
&+\binom{r}{\frac{r}{2}}\frac{A^{\frac{r}{2}}
(-B)^{\frac{r}{2}}}{1-(-1)^{\frac{r}{2}} x}\\
=&
\sum_{k=0}^{\frac{r}{2}-1}(-1)^k \binom{r}{k}
\frac{A^k B^{r-k}+A^{r-k}B^k-(-1)^k b^k(A^k B^{r-k}\al^{r-2k}+
A^{r-k} B^k\be^{r-2k})x}{1- (-1)^k b^kV_{r-2k}x+x^2}\\
&+\binom{r}{\frac{r}{2}}\frac{A^{\frac{r}{2}}
(-B)^{\frac{r}{2}}}{1-(-1)^{\frac{r}{2}} x}.
\end{split}
\]
\end{proof}
We can derive the following beautiful identities
\begin{cor}
\label{corU}
If $U_0=0$, then $A=B=\frac{U_1}{\al-\be}$ and
\[
\begin{split}
U(r,x)&=A^{r-1}
\sum_{k=0}^{\frac{r-1}{2}} \binom{r}{k}
\frac{b^k U_{r-2k} x}{1- (-b)^k V_{r-2k}x-x^2},\
\text{if $r$ odd}\\
U(r,x)&=A^r\sum_{k=0}^{\frac{r}{2}-1} (-1)^k \binom{r}{k} \frac{2-(-b)^k
V_{r-2k}x}{1-(-b)^k V_{r-2k}x+x^2}+
\binom{r}{\frac{r}{2}}
\frac{ (-1)^{\frac{r}{2}} A^r}{1-(-1)^{\frac{r}{2}} x},\
\text{if $r$ even.}
\end{split}
\]
\end{cor}
\begin{cor}
If $\{U_n\}_n$ is a non-degenerate
second-order recurrence sequence and $U_0=0$, then
\begin{eqnarray}
U(1,x)&=&  \frac{A^2U_1 x}{1-V_1x-x^2}\\
U(2,x)&=& \frac{-A^2(V_2+2)x(x-1)}{(x+1)(x^2-V_2 x+1)}\\
U(3,x) &=& \frac{A^4 U_1 x
\left( (a^2+2b)-2a^2 b x-(a^2+2b) x^2
\right)}{(1-V_3 x-x^2)(1+bV_1 x-x^2)}
\end{eqnarray}
\end{cor}
\begin{proof}
We use Corollary \refth{corU}.
The first two identities are straightforward.
Now,
\begin{eqnarray*}
U(3,x)&=&
A^4\left(\frac{U_3 x}{1-V_3 x-x^2}+\frac{bU_1 x}{1+bV_1 x-x^2} \right)\\
&=& A^4x\frac{U_3+bU_1+b(U_3V_1-U_1V_3)x-
(U_3+bU_1)x^2}{(1-V_3 x-x^2)(1+bV_1 x-x^2)}\\
&=&  \frac{A^4 U_1 x
\left( (a^2+2b)-2a^2 b x-(a^2+2b) x^2
\right)}{(1-V_3 x-x^2)(1+bV_1 x-x^2)},
\end{eqnarray*}
since $U_3+bU_1=(a^2+2b)U_1$ and $U_3V_1-U_1V_3=-2a^2 U_1$.
\end{proof}
\begin{rem}
If $U_n$ is the Fibonacci sequence, then $a=b=1$, and if $U_n$ is the
Pell sequence, then $a=2,b=1$.
\end{rem}

\section{Horadam's Theorem}

Horadam \cite{Horadam1994} found some closed forms for partial sums
$S_n=\displaystyle \sum_{i=1}^n P_i $,
$S_{-n}=\displaystyle \sum_{i=1}^n P_{-i}$,
where $P_n$ is
 the generalized Pell sequence, $P_{n+1}=2P_n+P_{n-1},P_1=p,P_2=q$.
 Let $p_n$ be the ordinary Pell sequence, with $p=1,q=2$, and
  $q_n$ be the sequence satisfying the same recurrence, with $p=1,q=3$.
 He proved
 \begin{thm}[Horadam]
 \label{thm_horadam}
 For any $n$,
\begin{eqnarray*}
&S_{4n}=q_{2n}(p q_{2n-1}+q q_{2n})+p-q;
&S_{4n-2}=q_{2n-1}(p q_{2n-2}+q q_{2n-1})\\
&S_{4n+1}= q_{2n}(p q_{2n}+q q_{2n+1})-q;
&S_{4n-1}= q_{2n}(p q_{2n-2}+q q_{2n-1})-q\\
&S_{-4n}=q_{2n}(-p q_{2n+2}+q q_{2n+1})+3p -q;
&S_{-4n+2}=q_{2n}(-p q_{2n}+qq_{2n-1})+2p\\
&S_{-4n+1}= q_{2n}(p q_{2n+1}-qq_{2n})+p;
&S_{-4n-1}= q_{2n+1}(p q_{2n+2}-qq_{2n+1})+2p-q.
\end{eqnarray*}
 \end{thm}

We observe that Horadam's theorem is a particular case of the partial sum
for a non-degenerate
 second-order recurrence sequence $U_n$.
 In fact,
  we find $\displaystyle S_{n,r}^U(x)=\sum_{i=0}^n U_i^r x^i$.
  For simplicity, we let $U_0=0$. Thus,
  $U_n=A(\al^n-\be^n)$ and $V_n=\al^n+\be^n$.
  We prove
  \begin{thm} We have
  \label{partial_thm}
 \begin{eqnarray}
S_{n,r}^U (x)= A^{r-1} \sum_{k=0}^{\frac{r-1}{2}} \binom{r}{k}
\frac{U_{r-2k} x-(-1)^{kn} U_{(r-2k)(n+1)} x^{n+1}-
(-1)^{k(n+1)} U_{(r-2k)n} x^{n+2}}{1-(-1)^k V_{r-2k} x-x^2}
 \end{eqnarray}
 if $r$ odd, and
 \begin{equation}
 \begin{split}
S_{n,r}^U(x)&= A^r \sum_{k=0}^{\frac{r-1}{2}} \binom{r}{k}
\frac{V_{r-2k} x-(-1)^{kn} V_{(r-2k)(n+1)} x^{n+1}-
(-1)^{k(n+1)} V_{(r-2k)n} x^{n+2}}{1-(-1)^k V_{r-2k} x+x^2}\\
&\hspace{2cm} +A^r \binom{r}{\frac{r}{2}} \frac{
(-1)^{\frac{r}{2} (n+1)} x^{n+1}-1}{(-1)^{\frac{r}{2}} x-1}
\end{split}
 \end{equation}
 if $r$ even.
  \end{thm}
\begin{proof}
We evaluate
\begin{eqnarray*}
S_{n,r}^U(x)&=&\sum_{i=0}^n \sum_{k=0}^r \binom{r}{k}
(A\al)^k (-A\be)^{r-k} x^i\\
&=& A^r \sum_{k=0}^r (-1)^{r-k} \binom{r}{k}
\sum_{i=0}^n (\al^k\be^{r-k} x)^i\\
&=& A^r \sum_{k=0}^r (-1)^{r-k}\binom{r}{k}
\frac{(\al^k\be^{r-k}x)^{n+1}-1}{\al^k\be^{r-k}x-1}.
 \end{eqnarray*}
 Assume $r$ odd. Then, associating $k\leftrightarrow r-k$, we get
\begin{eqnarray*}
S_{n,r}^U(x)&=&  A^r \sum_{k=0}^{\frac{r-1}{2}} (-1)^k \binom{r}{k}
\left(
\frac{(\al^{r-k}\be^k x)^{n+1}-1}{\al^{r-k}\be^k x-1}-
\frac{(\al^k\be^{r-k}x)^{n+1}-1}{\al^k\be^{r-k}x-1}
\right)\\
&=&  A^r \sum_{k=0}^{\frac{r-1}{2}} (-1)^k \binom{r}{k}
\frac{(\al^k\be^{r-k}x-1)(\al^{(r-k)(n+1)}\be^{k(n+1)}x^{n+1}-1)}{}\\
&&\hspace{2cm}\frac{-
(\al^{r-k}\be^k x-1)(\al^{k(n+1)}
\be^{(r-k)(n+1)}x^{n+1}-1)}{(\al^k\be^{r-k}x-1)(\al^{r-k}\be^k x-1)}
\\
&=&  A^r \sum_{k=0}^{\frac{r-1}{2}} (-1)^k \binom{r}{k}
\frac{\al^{k(n+1)-kn}\be^{r+kn}x^{n+2}}{}\\
&&\hspace{2.2cm}\frac{
-\al^{(r-k)(n+1)}\be^{k(n+1)}x^{n+1}-\al^k\be^{r-k}x}{}\\
&&\hspace{2.4cm}\frac{
-\al^{r+kn}\be^{r(n+1)-kn}x^{n+2}+\al^{r-k}\be^k x}{}\\
&&\hspace{2.6cm}\frac{
-\al^{k(n+1)}\be^{(r-k)(n+1)}\be^{(r-k)(n+1)} x^{n+1}}
{1-(-1)^k(\al^{r-2k}+\be^{r-2k})+\al^r\be^r x^2}
 \end{eqnarray*}
\begin{eqnarray*}
&=&  A^r \sum_{k=0}^{\frac{r-1}{2}} (-1)^k \binom{r}{k}
\frac{(-1)^k(\al^{r-2k}-\be^{r-2k}) x-(-1)^{k(n+1)}
(\al^{(r-2k)(n+1)}}{}\\
&&\hspace{2cm}\frac{
 -\be^{(r-2k)(n+1)}) x^{n+1}+
(-1)^{r+kn)}(\al^{(r-2k)n} -
\be^{(r-2k)n})x^{n+2}}{1-(-1)^k V_{r-2k} x-x^2}\\
&=&  A^{r-1} \sum_{k=0}^{\frac{r-1}{2}} \binom{r}{k}
\frac{U_{r-2k} x-(-1)^{kn} U_{(r-2k)(n+1)} x^{n+1}-
(-1)^{k(n+1)} U_{(r-2k)n} x^{n+2}}{1-(-1)^k V_{r-2k} x-x^2}.
 \end{eqnarray*}

Assume $r$ even. Then, as before, associating $k\leftrightarrow r-k$,
except for the middle term, we get
\begin{eqnarray*}
S_{n,r}^U(x)&=&  A^r \sum_{k=0}^{\frac{r-1}{2}} (-1)^k \binom{r}{k}
\frac{(-1)^k(\al^{r-2k}+\be^{r-2k}) x-(-1)^{k(n+1)}
(\al^{(r-2k)(n+1)}}{}\\
&&\hspace{2cm}\frac{
 +\be^{(r-2k)(n+1)}) x^{n+1}+
(-1)^{r+kn)}(\al^{(r-2k)n} +
\be^{(r-2k)n})x^{n+2}}{1-(-1)^k V_{r-2k} x+x^2}\\
&&\hspace{2cm} +A^r \binom{r}{\frac{r}{2}} \frac{
(-1)^{\frac{r}{2} (n+1)} x^{n+1}-1}{(-1)^{\frac{r}{2}} x-1}
\\
 &=& A^r \sum_{k=0}^{\frac{r-1}{2}} \binom{r}{k}
\frac{V_{r-2k} x-(-1)^{kn} V_{(r-2k)(n+1)} x^{n+1}-
(-1)^{k(n+1)} V_{(r-2k)n} x^{n+2}}{1-(-1)^k V_{r-2k} x+x^2}\\
&&\hspace{2cm} +A^r \binom{r}{\frac{r}{2}} \frac{
(-1)^{\frac{r}{2} (n+1)} x^{n+1}-1}{(-1)^{\frac{r}{2}} x-1}.
\end{eqnarray*}
\end{proof}
Taking $r=1$, we get the partial sum for any non-degenerate second-order
recurrence sequence, with $U_0=0$,
\begin{cor}
$\displaystyle S_{n,1}^U(x)=
\frac{x\left( U_1-U_{n+1} x^n-U_n x^{n+2} \right)}{1-V_1 x-x^2}$
\end{cor}

\begin{rem}
Horadam's theorem follows easily, since
$S_n=S_{n,1}^{P}(1)$.
Also $S_{-n}$ can be found without difficulty, by observing
that $P_{-n}=p p_{-n-2}+q p_{-n-1}=
-p (-1)^{n+2} p_{n+2}-q(-1)^{n+1}p_{n+1}$, and using $S_{n,1}^{p}(-1)$.
\end{rem}

\section{Weighted Combinatorial Sums}

In \cite{Vajda1989} there are quite a few identities of the form
\(\displaystyle \sum_{i=0}^n\binom{n}{i} F_i=F_{2n}\), or
\(\displaystyle \sum_{i=0}^n\binom{n}{i} F_i^2\), which
is $5^{[\frac{n-1}{2}]} L_n$ if
$n$ even, and $5^{[\frac{n-1}{2}]} F_n$, if $n$ odd.
 A natural question is: {\em for fixed $r$, what is the closed form for
 the weighted sum
$\sum_{i=0}^n\binom{n}{i} F_i^r$ (if it exists)?} We are able to answer the
previous question, not only for the Fibonacci sequence, but also
for any second-order recurrence sequences. Let
$\displaystyle S_{r,n}(x)=\sum_{i=0}^n\binom{n}{i} U_i^r x^i$.
\begin{thm}
\label{thmS}
We have
\[
S_{r,n}(x)=
 \sum_{k=0}^r \binom{r}{k} A^k (-B)^{r-k} (1+\al^k \be^{r-k}x)^n.
\]
Moreover, if $U_0=0$, then
$\displaystyle S_{r,n}(x)=
A^r\sum_{k=0}^r (-1)^{r-k}\binom{r}{k} (1+\al^k\be^{r-k}x)^n$.
\end{thm}
\begin{proof}
Let
\begin{eqnarray*}
S_{r,n}(x)&=&
\sum_{i=0}^n\binom{n}{i} \sum_{k=0}^r \binom{r}{k}
(A\al^i)^k(-B\be^i)^{r-k}x^i\\
&=&
\sum_{k=0}^r \binom{r}{k} A^k (-B)^{r-k} \sum_{i=0}^n\binom{n}{i}
(\al^k\be^{r-k}x)^i\\
&=& \sum_{k=0}^r \binom{r}{k} A^k (-B)^{r-k} (1+\al^k \be^{r-k}x)^n
\end{eqnarray*}
If $U_0=0$, then $A=B$, and
\(
S_{r,n}(x)=A^r \sum_{k=0}^r (-1)^{r-k}\binom{r}{k} (1+\al^k\be^{r-k}x)^n
\)
\end{proof}


Studying Theorem \refth{thmS},
we observe that we get nice sums involving the
Fibonacci and Lucas sequences (or any such sequence, for that matter),
 if we
are able to express 1 plus/minus a power of $\al,\be$ as the same
multiple of a power of $\al$, respectively $\be$.
The following lemma turns out to be very useful.
\begin{lem}
\label{th_dodd}
The following identities are true
\begin{equation}
\label{eq_dodd}
\begin{split}
\al^{2s}-(-1)^s=&\sqrt{5} \al^s F_s\\
\be^{2s}-(-1)^s=&-\sqrt{5} \be^s F_s
\\
\al^{2s}+(-1)^s=& L_s\al^s \\
\be^{2s}+(-1)^s=& L_s\be^s.
\end{split}
\end{equation}
\end{lem}
\begin{proof}
Straightforward using the Binet formula for $F_s$ and $L_s$.
\end{proof}

\begin{thm}
We have
\begin{eqnarray}
S_{4r,n}(1)&=& 5^{-2r}\left(
\sum_{k=0}^{2r-1}(-1)^{k(n+1)}\binom{4r}{k}
 L_{2r-k}^n L_{(2r-k)n}+ \binom{4r}{2r} 2^n\right)\\
S_{4r+2,n}(1)&=&5^{\frac{n+1}{2}-(2r+1)} \sum_{k=0}^{2r}\binom{4r+2}{k}
F_{2r+1-k}^n F_{n(2r+1-k)},\ \text{if $n$ odd}\\
S_{4r+2,n}(1)&=& 5^{\frac{n}{2}-(2r+1)}
\sum_{k=0}^{2r}(-1)^k \binom{4r+2}{k}
F_{2r+1-k}^n L_{n(2r+1-k)}\ \text{if $n$ even}.
\end{eqnarray}
\end{thm}
\begin{proof}
We use Theorem \refth{thmS}. Associating $k\leftrightarrow 4r+2-k$,
except for the middle term in $S_{4r+2,n}(1)$, we obtain
\begin{equation}
\label{eq_s9}
\begin{split}
&S_{4r+2,n}(1)=5^{-(2r+1)}
\sum_{k=0}^{2r}(-1)^k\binom{4r+2}{k}
\left( (1+\al^k\be^{4r+2-k})^n+ (1+\al^{4r+2-k}\be^k)^n \right)\\
&=5^{-(2r+1)}
\sum_{k=0}^{2r}(-1)^k\binom{4r+2}{k}
\left( (1+(-1)^k\be^{4r+2-2k})^n+ (1+(-1)^k\al^{4r+2-2k})^n \right)\\
&=5^{-(2r+1)}
\sum_{k=0}^{2r}(-1)^{k(n+1)}\binom{4r+2}{k}
\left( ((-1)^k+\be^{2(2r+1-k)})^n+ ((-1)^k+\al^{2(2r+1-k)})^n \right).
\end{split}
\end{equation}
We did not insert the middle term, since it is equal to
\[
\begin{split}
&5^{-(2r+1)} (-1)^{2r+1}\binom{4r+2}{2r+1}(1+\al^{2r+1}\be^{2r+1})^n\\
&=5^{-(2r+1)} (-1)^{2r+1}\binom{4r+2}{2r+1}(1+(-1)^{2r+1})^n=0.
\end{split}
\]

Assume first that $n$ is odd.
Using \refeq{eq_dodd} into \refeq{eq_s9}, and observing that
$\al^{2(2r+1-k)}-(-1)^{2r+1-k}=\al^{2(2r+1-k)}+(-1)^k$,
we get
\[
\begin{split}
S_{4r+2,n}(1)=5^{-(2r+1)}
&\sum_{k=0}^{2r}(-1)^{(n+1)k}\binom{4r+2}{k}
5^{\frac{n}{2}}F_{2r+1-k}^n\\
&
\left(
(-1)^n \be^{n(2r+1-k)}+\al^{n(2r+1-k)}
\right)\\
=5^{-(2r+1)}
&\sum_{k=0}^{2r}\binom{4r+2}{k}
5^{\frac{n+1}{2}}F_{2r+1-k}^n F_{n(2r+1-k)}
\end{split}
\]

Assume $n$ even. As before,
\[
\begin{split}
S_{4r+2,n}(1)=5^{-(2r+1)}
&\sum_{k=0}^{2r}(-1)^{(n+1)k}\binom{4r+2}{k}
5^{\frac{n}{2}}F_{2r+1-k}^n\\
&
\left(
(-1)^n \be^{n(2r+1-k)}+\al^{n(2r+1-k)}
\right)\\
=5^{-(2r+1)}
&\sum_{k=0}^{2r}(-1)^k \binom{4r+2}{k}
5^{\frac{n}{2}}F_{2r+1-k}^n L_{n(2r+1-k)}
\end{split}
\]
In the same way, associating $k\leftrightarrow 4r-k$,
except for the middle term,
\begin{equation}
\begin{split}
S_{4r,n}(1)&=5^{-2r}
\sum_{k=0}^{2r-1}(-1)^k\binom{4r}{k}
\left( (1+\al^k\be^{4r-k})^n+ (1+\al^{4r-k}\be^k)^n \right)+
5^{-2r}\binom{4r}{2r} 2^n\\
&=5^{-2r}
\sum_{k=0}^{2r-1}(-1)^{k(n+1)}\binom{4r}{k}
\left( \left((-1)^k+\be^{2(2r-k)}\right)^n+
\left((-1)^k+\al^{2(2r-k)}\right)^n \right)\\
&\qquad \qquad\qquad + 5^{-2r}\binom{4r}{2r} 2^n\\
&=5^{-2r}\left(
\sum_{k=0}^{2r-1}(-1)^{k(n+1)}\binom{4r}{k}
\left( L_{2r-k}^n\be^{(2r-k)n}+L_{2r-k}^n\al^{(2r-k)n} \right)
+\binom{4r}{2r} 2^n\right)\\
&=5^{-2r}\left(
\sum_{k=0}^{2r-1}(-1)^{k(n+1)}\binom{4r}{k}
 L_{2r-k}^n L_{(2r-k)n}+ \binom{4r}{2r} 2^n\right).
\end{split}
\end{equation}
\end{proof}
\begin{rem}
In the same manner we can find
$\displaystyle \sum_{i=0}^n \binom{n}{i} U_{pi}^r x^i$.
\end{rem}
As a consequence of the previous theorem,
for the even cases, and working out
the details for the odd cases we get
\begin{cor} We have
\begin{eqnarray*}
&& \sum_{k=0}^{n} \binom{n}{i} F_i =F_{2n}\\
&& \sum_{k=0}^{2n} \binom{2n}{i} F_i^2 =5^{n-1} L_{2n}\\
&& \sum_{k=0}^{2n+1} \binom{2n+1}{i} F_i^2 =5^n F_{2n+1}
\\
&& \sum_{k=0}^{n} \binom{n}{i} F_i^3 =\frac{1}{5}(2^n F_{2n}+3F_n)
\\
&& \sum_{k=0}^{n} \binom{n}{i} F_i^4 =\frac{1}{25}
(3^nL_{2n}-4(-1)^nL_n+6\cdot 2^n).
\end{eqnarray*}
\end{cor}
\begin{proof}
The second, third and fifth identities follow from the previous theorem. Now,
 using Theorem \refth{thmS}, with $A=\frac{1}{\sqrt{5}}$,
 we get
\begin{eqnarray*}
S_{1,n}(1)&=&\frac{1}{\sqrt{5}}
\sum_{k=0}^1 (-1)^{1-k}\binom{1}{k} (1+\al^k\be^{1-k})^n\\
&=& \frac{1}{\sqrt{5}}(-(1+\be)^n+(1+\al)^n )=
\frac{1}{\sqrt{5}}(\al^{2n}-\be^{2n})=F_{2n}.
\end{eqnarray*}
Furthermore, the fourth identity follows from
\begin{eqnarray*}
S_{3,n}(1)&=&\frac{1}{5\sqrt{5}}
\sum_{k=0}^3 (-1)^{3-k}\binom{3}{k} (1+\al^k\be^{3-k})^{n}\\
&=& \frac{1}{5\sqrt{5}}(-(1+\be^3)^{n}+3(1+\al\be^2)^{n}-
3(1+\al^2\be)^{n}+(1+\al^3)^{n})\\
&=& \frac{1}{5\sqrt{5}}(-(2\be^2)^{n}+3\al^{n}-3\be^{n}+(2\al^2)^{n})\\
&=& \frac{1}{5}(2^n F_{2n}+3 F_n),
\end{eqnarray*}
since $1+\be^3=2\be^2,\, 1+\al^3=2\al^2$.
\end{proof}


We remark the following
\begin{cor}
We have, for any $n$,
\item[(i)]  $2^n F_{2n}+3F_n\equiv 0\pmod 5$
\item[(ii)] $3^nL_{2n}-4(-1)^nL_n+6\cdot 2^n\equiv 0\pmod{5^2}$
\item[(iii)] $\displaystyle \sum_{k=0}^{2r}\binom{4r+2}{k}
F_{2r+1-k}^n F_{n(2r+1-k)}\equiv 0
\pmod{5^{4r+2-\frac{n-1}{2}}}$,
if $n$ is odd, $n\leq 8r+3$.
\item[(iv)] $\displaystyle \sum_{k=0}^{2r}(-1)^k
 \binom{4r+2}{k}
F_{2r+1-k}^n L_{n(2r+1-k)}\equiv 0
\pmod{5^{4r+2-\frac{n}{2}}}$,
if $n$ is even, $n\leq 8r+2$.
\item[(v)]
$ \displaystyle
\sum_{k=0}^{2r-1}(-1)^{k(n+1)}\binom{4r}{k}
 L_{2r-k}^n L_{(2r-k)n}+ \binom{4r}{2r} 2^n\equiv 0\pmod{5^{2r}}$.
\end{cor}

Taking other values for $x$ (as desired) in Theorem \refth{thmS},
for instance, $x=-1$
 and working out the details, we get the following
\begin{thm}
\label{thm_bin2}
We have
\begin{eqnarray*}
S_{4r,n}(-1)&=&5^{\frac{n}{2}-2r}
\sum_{k=0}^{2r-1}(-1)^{k}F_{2r-k}^n L_{(4r-2k)n}\binom{4r}{k},
\ \text{if $n$ even}\\
S_{4r,n}(-1)&=&-5^{\frac{n+1}{2}-2r}
\sum_{k=0}^{2r-1}F_{2r-k}^n F_{(4r-2k)n}\binom{4r}{k},
\ \text{if $n$ odd}\\
S_{4r+2,n}(-1)&=&5^{-(2r+1)}\left(
\sum_{k=0}^{2r}(-1)^{k(n+1)+n}\binom{4r+2}{k}
L_{2r+1-k}^n L_{(2r+1-k)n} -  2^n \binom{4r+2}{2r+1}\right).
\end{eqnarray*}
\end{thm}
\begin{proof}
We use $x=-1$ in Theorem \refth{thmS}. Associating
$k\leftrightarrow 4r+2-k$ in $S_{4r+2,n}(-1)$, we obtain
\begin{equation*}
\begin{split}
&S_{4r+2,n}(-1)=5^{-(2r+1)}
\sum_{k=0}^{2r}(-1)^k\binom{4r+2}{k}
\left( (1-\al^k\be^{4r+2-k})^n+ (1-\al^{4r+2-k}\be^k)^n \right)\\
&\hspace{6cm} - 5^{-(2r+1)} 2^n \binom{4r+2}{2r+1}\\
&=5^{-(2r+1)}
\sum_{k=0}^{2r}(-1)^k\binom{4r+2}{k}
\left( (1-(-1)^k\be^{4r+2-2k})^n+ (1-(-1)^k\al^{4r+2-2k})^n \right)\\
&\hspace{6cm} - 5^{-(2r+1)} 2^n \binom{4r+2}{2r+1}
\end{split}
\end{equation*}
\begin{equation*}
\begin{split}
&=5^{-(2r+1)}
\sum_{k=0}^{2r}(-1)^{k(n+1)}\binom{4r+2}{k}
\left( ((-1)^k-\be^{2(2r+1-k)})^n+ ((-1)^k-\al^{2(2r+1-k)})^n \right)\\
&\hspace{6cm} - 5^{-(2r+1)} 2^n \binom{4r+2}{2r+1}\\
&=5^{-(2r+1)}
\sum_{k=0}^{2r}(-1)^{k(n+1)+n}\binom{4r+2}{k}
L_{2r+1-k}^n L_{(2r+1-k)n}
 - 5^{-(2r+1)} 2^n \binom{4r+2}{2r+1},
\end{split}
\end{equation*}
since $(-1)^k-\be^{4r+2-2k}=-L_{2r+1-k}\be^{2r+1-k}$ and
$(-1)^k-\al^{4r+2-2k}=-L_{2r+1-k}\al^{2r+1-k}$, by Lemma~ \refth{th_dodd}.
In the same way, associating $k\leftrightarrow 4r-k$, with the
middle term zero,
\begin{equation*}
\begin{split}
S_{4r,n}(1)&=5^{-2r}
\sum_{k=0}^{2r-1}(-1)^k\binom{4r}{k}
\left( (1-\al^k\be^{4r-k})^n+ (1-\al^{4r-k}\be^k)^n \right)\\
&=5^{-2r}
\sum_{k=0}^{2r-1}(-1)^{k(n+1)}\binom{4r}{k}
\left( \left((-1)^k-\be^{2(2r-k)}\right)^n+
\left((-1)^k-\al^{2(2r-k)}\right)^n \right)
\end{split}
\end{equation*}
\begin{equation*}
\begin{split}
&=5^{-2r}
\sum_{k=0}^{2r-1}(-1)^{k(n+1)}\binom{4r}{k}
\left(5^{\frac{n}{2}} F_{2r-k}^n\be^{(2r-k)n}+
5^{\frac{n}{2}}(-1)^n F_{2r-k}^n\al^{(2r-k)n} \right)\\
&=5^{\frac{n}{2}-2r}
\sum_{k=0}^{2r-1}(-1)^{k(n+1)+n}F_{2r-k}^n\binom{4r}{k}
(\al^{(2r-k)n}+(-1)^n\be^{(2r-k)n}),
\end{split}
\end{equation*}
since  $(-1)^k-\be^{4r-2k}=\sqrt{5}F_{2r-k}\be^{2r-k}$ and
$(-1)^k-\al^{4r-2k}=-\sqrt{5}F_{2r-k}\al^{2r-k}$,
by Lemma~ \refth{th_dodd}. Therefore, for $n$ even,
\(\displaystyle
S_{4r,n}(1)=5^{\frac{n}{2}-2r}
\sum_{k=0}^{2r-1}(-1)^{k(n+1)+n}F_{2r-k}^n L_{(4r-2k)n}\binom{4r}{k},
\)
and for $n$ odd,
\(\displaystyle
S_{4r,n}(1) =5^{\frac{n+1}{2}-2r}
\sum_{k=0}^{2r-1}(-1)^{k(n+1)+n}F_{2r-k}^n F_{(4r-2k)n}\binom{4r}{k}.
\)
\end{proof}
A consequence for even powers and a similar idea for odd powers produces
\begin{cor}
We have
\begin{eqnarray*}
\sum_{i=0}^n (-1)^i\binom{n}{i} F_i &=& -F_n\\
\sum_{i=0}^n (-1)^i\binom{n}{i} F_i^2&=&
\frac{1}{5}\left((-1)^n L_n-2^{n+1}\right)\\
\sum_{i=0}^n (-1)^i\binom{n}{i} F_i^3&=& \frac{1}{5}
\left((-2)^n F_n-3 F_{2n}\right)
\\
\sum_{i=0}^{n} (-1)^i\binom{n}{i} F_i^4&=&
5^{\frac{n}{2}-2}(L_{2n}-L_{n}),\ \text{if $n$ even}\\
\sum_{i=0}^{n} (-1)^i\binom{n}{i} F_i^4&=&
-5^{\frac{n+1}{2}-2} (F_{2n}+4 F_{n}), \ \text{if $n$ odd}.
\end{eqnarray*}
\end{cor}
\begin{proof}
The first identity is simple application of Theorem \refth{thmS}.
The identities for even powers are consequences of Theorem \refth{thm_bin2}.
Now, using Theorem~ \refth{thmS}, we get
\[
\begin{split}
S_{3,n}(-1)&=\frac{1}{5\sqrt{5}}\left(-(1-\be^3)^n+
3(1-\al\be^2)^n-3(1-\al^2\be)^n+(1-\al^3)^n \right)\\
&=\frac{1}{5\sqrt{5}}\left(
(-2)^n \be^n+3\be^{2n}-3\al^{2n} +(-2)^n \al^n
\right)=\frac{1}{5}
((-2)^n F_n-3 F_{2n}),
\end{split}
\]
since $1-\be^3=-2\be,\, 1-\al^3=-2\al$.
\end{proof}
We remark the following
\begin{cor}
We have, for any $n$,
$(-1)^n L_n-2^{n+1}\equiv 0\pmod{5}$ and
$(-2)^n F_n-3 F_{2n}\equiv 0\pmod 5.$
\end{cor}

\noindent AMS Classification Numbers: 05A10, 05A19, 11B37, 11B39, 11B65,
11B83

\end{document}